\documentclass[12pt]{article}

\newif\ifpdf
\ifx\pdfoutput\undefined
    \pdffalse       
\else
    \pdfoutput=1    
    \pdftrue
\fi

\usepackage{amsmath,amssymb,amsthm}
\usepackage[numeric]{amsrefs} 

\ifpdf
    \usepackage[pdftex]{graphicx}
    \usepackage[backref, colorlinks=true, pdfstartview=FitH, linkcolor=blue,
        citecolor=blue, urlcolor=blue]{hyperref}
\else
    \usepackage{graphicx}
    \usepackage{hyperref}
\fi

\title{Fraenkel's Partition and Brown's Decomposition}
\author{Kevin O'Bryant \\ University of California - San Diego \\
        \href{mailto:kobryant@math.ucsd.edu}{kobryant@math.ucsd.edu} \\
        \href{http://www.math.ucsd.edu/~kobryant}{www.math.ucsd.edu/$\sim$kobryant}}
\date{\today}

    \newcommand{\N}{{\mathbb N}}
    \newcommand{\Z}{{\mathbb Z}}
    \newcommand{\Q}{{\mathbb Q}}
    \newcommand{\R}{{\mathbb R}}
    \newcommand{\C}{{\mathbb C}}

    \newcommand{\B}{{\cal B}}

    \newcommand{\floor}[1]{\left\lfloor #1 \right\rfloor}
    \newcommand{\ceiling}[1]{\left\lceil #1 \right\rceil}

    \newcommand{\fp}[1]{\left\{ #1 \right\}}
    \newcommand{\MR}[1]{\href{http://www.ams.org/mathscinet-getitem?mr=#1}{\mbox{\bf MR~#1}}}

    \newcommand{\arc}[1]{\overline{#1}}

    \newtheorem{thm}{Theorem}
    \newtheorem{lem}[thm]{Lemma}

    \newenvironment{proofof}[1]{\medskip\noindent{\em Proof of #1.}}{\qed\medskip}
    \newenvironment{namedtheorem}[1]{\medskip \noindent {\bf #1:}\begin{em}}{\end{em}\medskip}

    \newcommand{\al}{\alpha}
    \newcommand{\alpr}{\alpha^{\prime}}
    \newcommand{\alprpr}{\frac{\ceiling{q\alpr}}{q}}
    \newcommand{\ap}{a^{\prime}}
    \newcommand{\be}{\beta}
    \newcommand{\bepr}{\beta^{\prime}}
    
    \newcommand{\bp}{b^{\prime}}

\begin{document}
    \ifpdf
       \DeclareGraphicsExtensions{.pdf,.jpg,.mps,.png}
    \fi
\maketitle

\begin{abstract}
We give short proofs of Fraenkel's Partition Theorem and Brown's Decomposition. Denote the sequence
$\left(\floor{(n-\alpr)/\al}\right)_{n=1}^\infty$ by $\B(\al,\alpr)$, a so-called Beatty sequence. Fraenkel's
Partition Theorem gives necessary and sufficient conditions for $\B(\al,\alpr)$ and $\B(\be,\bepr)$ to tile
the positive integers, i.e., for $\B(\al,\alpr) \cap \B(\be,\bepr)=\emptyset$ and $\B(\al,\alpr) \cup
\B(\be,\bepr)=\N$. Fix $\al\in(0,1)$, and let $c_k=1$ if $k\in \B(\al,0)$, and $c_k=0$ otherwise, i.e.,
$c_k=\floor{(k+1)\al}-\floor{k\al}$. For a positive integer $m$ let $C_m$ be the binary word $c_1c_2c_3\cdots
c_m$. Brown's Decomposition gives integers $q_1,q_2,\dots$, {\em independent of $m$ and growing at least
exponentially}, and integers $t,z_0,z_1,z_2,\dots,z_t$ (depending on $m$) such that
$C_m=C_{q_t}^{z_t}C_{q_{t-1}}^{z_{t-1}} \cdots C_{q_1}^{z_1}C_{q_0}^{z_0}$. In other words, Brown's
Decomposition gives a sparse set of initial segments of $C_\infty$ and an explicit decomposition of $C_m$
(for every $m$) into a product of these initial segments.
\end{abstract}

 \thispagestyle{empty}

 \pagebreak

 \tableofcontents

 \pagebreak

\section{Introduction}

In 1894, Rayleigh \cite{Rayleigh} observed that
\begin{quote}
``If $x$ be an incommensurable number less than unity, one of the series of quantities $m/x, m/(1-x)$,
where $m$ is a whole number, can be found which shall lie between any given consecutive integers, and
but one such quantity can be found.''
\end{quote}
S. Beatty \cite{1926.Beatty} posed this as a Monthly problem in 1926, and it has come to be known as Beatty's
Theorem.

The Beatty sequence with density $\al$ and offset $\alpr$ is defined by
    \begin{equation}\label{BeattyDefinition}
        \B(\al,\alpr) := \left( \floor{\frac{n-\alpr}{\al}} \right)_{n=1}^\infty,
    \end{equation}
where $\floor{x}$ is the floor of $x$. When the second argument is 0 we omit it from our notation, i.e.,
$\B(\al):=\B(\al,0)$. We write $\fp{x}:=x-\floor{x}$ for the fractional part of $x$, and $\ceiling{x}$ for
the smallest integer larger than or equal to $x$. We say that two sequences {\em tile} a set $S$ if they are
disjoint and their union is $S$.


For example, we now state Beatty's Theorem in this language.

\begin{namedtheorem}{Beatty's Theorem}
The sequences $\B(\al)$ and $\B(\be)$ tile $\N$ if and only if $0<\al<1$, $\al+\be=1$, and $\al$ is
irrational.
\end{namedtheorem}

Beatty sequences arise in a number of areas, including Computer Graphics, Signal Processing, Automata,
Quasicrystals, Combinatorial Games, and Diophantine Approximation. They are natural counterparts to
Kronecker's fractional part sequences. There is the obvious connection of $\B(\al,\alpr)$ to
$\left(\fp{\frac{n-\alpr}{\al}}\right)_{n=1}^\infty$, but also a more subtle connection to
$\left(\fp{n\al+\alpr}\right)_{n=1}^\infty$ which we exploit here to give simple proofs of Fraenkel's
Partition and Brown's Decomposition.

Beatty sequences generalize arithmetic progressions, which correspond to the special case $\al^{-1}\in\Z$.
Most work on Beatty sequences has the aim of extending some known result on APs. For example, the Chinese
Remainder Theorem identifies precisely when APs are disjoint; Fraenkel's Partition (stated precisely in
Section~\ref{Sec.Statement.Fraenkel}) identifies precisely when two Beatty sequences are disjoint and have
union $\N$. The situation for more than 2 sequences is inadequately understood; see
\cites{MR2001f:11039,MR2001h:11011} for an up-to-date survey of knowledge in that direction.

Fraenkel proved his elegant generalization of Beatty's Theorem in 1969. Although his argument is well
motivated and geometric, it is rather long and hampered by unfortunate notation. Skolem~\cite{MR19:1159i}
attempted to deal with the special case of Beatty sequences with irrational densities. Alas, both his
statement of the theorem and his proof are incorrect (this is discussed in \cite{Fraenkel}). Borwein \&
Borwein~\cite{MR94i:11044} give ``a new proof of a theorem of Fraenkel''. They write, ``Our proof, once we
have developed the other machinery of this paper, is considerably shorter.'' Their proof is indeed short and
the ``other machinery'' consists only of two straightforward functional equations relating a pair of
generating functions. Unfortunately, the ``theorem of Fraenkel'' that they prove is a {\em very}\ special
case of Fraenkel's Partition. In Section~\ref{Sec.Statement.Fraenkel}, we give the statements put forward by
Skolem, Fraenkel, and Borwein \& Borwein.

Skolem's work is incorrect, Fraenkel's proof is long, and the Borwein brothers shied away from proving the
full theorem. For these reasons, the author feels that there is room in the literature for another proof of
Fraenkel's Partition, provided that it is correct, short, and complete. We state Fraenkel's Partition in
Section~\ref{Sec.Statement.Fraenkel}, define some relevant notation in Section~\ref{Sec.Preliminaries}, and
give the proof in Section~\ref{Sec.Proof.Fraenkel}. Because of the theorem's history of incorrect proofs and
inadequate statements, we have perhaps erred on the side of including too much detail. To ease the work of a
casual reader, the main ideas are given in Section~\ref{Sec.Special.Fraenkel}.

We now introduce Brown's Decomposition. Fix an $\al\in(0,1)$, and set
    $$
    c_k :=
            \left\{%
            \begin{array}{ll}
                1, & \hbox{$k \in \B\left(\al\right)$;} \\
                0, & \hbox{otherwise,} \\
            \end{array}%
            \right.
    $$
and let $C_m$ be the binary word $c_1c_2\dots c_m$. The word $C_\infty$ is called the characteristic word
with density $\al$. If $\al^{-1}=q\in\N$, then $c_k=1$ if and only if $k\equiv \floor{-q\alpr} \pmod{q}$. In
particular, $c_k=c_{q+k}$ for every $k$, whence $C_m = C_q^{\floor{m/q}}C_{m-q\floor{m/q}}$. In fact, if
$\al=\frac{a}{q}$, then the sequence $c_1,c_2,\dots$ is periodic with period $q$ and so
$C_m=C_q^{\floor{m/q}}C_{m-q\floor{m/q}}$.

If $\al$ is irrational, then the sequence $c_1,c_2,\dots$ is not periodic. But $\al$ is near to rationals,
and so an initial segment of the sequence will appear to be periodic. Brown's Decomposition (stated precisely
in Section~\ref{Sec.Statement.Brown}) is a quantitative description of this, given in terms of the
convergents of the continued fraction of $\al$.

We comment that this `almost-periodicity' is precisely what makes Beatty sequences interesting to
quasicrystallographers. If $\al\not\in\Q$, then $C_\infty$ is the prototypical example of a Sturmian Word.
Much of the literature on Beatty sequences is couched in the equivalent (and often more convenient) language
of Sturmian Words, especially the literature analyzing quasicrystallographic properties.

The first steps toward Brown's Decomposition were made in the oft-cited work of Stolarsky \cite{MR56:2908}.
He was studying functions $h$ for which
    $c_1 c_2 \cdots = h(c_1) h(c_2) \cdots.$
A nontrivial example is worth a thousand words: set $\al=\frac{\sqrt{5}-1}{2}$, in which case
    $$c_1 c_2 c_3 \cdots = {101}\,{101}\,{011}\,{011}\,{010} \cdots$$
is the ``Fibonacci word''. With $h$ defined by $h(1)=10, h(0)=1$ we have
\begin{align*}
 h(c_1)h(c_2)h(c_3) \cdots
    &= h(1)h(0)h(1) \cdots \\
    &= 10\,1\,10 \cdots \\
    &= c_1 c_2 c_3 c_4 c_5 \cdots.
\end{align*}
In the early 1990s, T. C. Brown found the useful and succinct decomposition of $C_m$ using morphisms. His
proof is nicely exposited in the new book of Allouche \& Shallit~\cite{Allouche.Shallit}. We remark that the
properties of morphisms received a great deal of attention throughout the 1990s. Excellent accounts of the
current theory are given in both \cite{Allouche.Shallit} and in the recent book of Lothaire \cite[Chapter
2]{Lothaire}.

Here, we give a short direct proof of Brown's Decomposition. Our proof relies on the same characterization of
$\B(\al)$ in terms of the fractional part sequence $\fp{k\al}$ that we use in our proof of Fraenkel's
Partition. We also need a well-known theorem from continued fractions (which is also used in Brown's proof).
The Decomposition is stated precisely in Section~\ref{Sec.Statement.Brown}, some notation is introduced in
Section~\ref{Sec.Preliminaries}, and the proof is given in Section~\ref{Sec.Proof.Brown}.

\section{Statement of Fraenkel's Partition}\label{Sec.Statement.Fraenkel}

\begin{namedtheorem}{Fraenkel's Partition Theorem}
The sequences $\B(\al,\alpr)$ and $\B(\be,\bepr)$ tile $\N$ if and only if the following five conditions are
satisfied.
\begin{enumerate}
 \item
    $0<\al<1$.
 \item
    $\al+\be=1$.
 \item
    $0\le\al+\alpr \le 1$.
 \item
    If $\al$ is irrational, then $\alpr+\bepr=0$ and $k\al+\alpr\not\in\Z$ for $2\le k\in \N$.
 \item
    If $\al$ is rational (say $q\in \N$ is minimal with $q\al \in \N$), then
    $\frac1q \le \al+\alpr$ and $\ceiling{q\alpr}+\ceiling{q\bepr}=1.$
\end{enumerate}
\end{namedtheorem}

We note first Conditions~{\em 1--5} are symmetric in $\al$ and $\be$. At first glance this is not the case;
for example, it is not clear that Conditions~{\em 1--5} imply that $0\le\be+\bepr\le1$. Our proof of
Fraenkel's Partition begins by proving the claimed symmetry.

Note also that if one divides the equation in Condition~{\em 5} by $q$, and then takes the limit as
$q\to\infty$, one obtains the equation $\alpr+\bepr=0$ of Condition~{\em 4}. This hints that the irrational
case can be derived as a limit of the rational case. However, the presence of the additional clause
``$k\al+\alpr\not\in\Z$'' of Condition~{\em 4} indicates that this approach is not trivial. The proof given
here handles the two cases separately. In Section~\ref{Sec.nonstandard}, we derive the irrational case from
the rational case using nonstandard analysis.

If one wishes to consider tilings of $\{N+1,N+2,\dots\}$ instead of $\{1,2,\dots\}$, it is not difficult to
adapt our statement. Indeed, $\B(\al,\alpr)$ and $\B(\be,\bepr)$ tile $\N$ if and only if
$\B(\al,\alpr-N\al)$ and $\B(\be,\bepr-N\be)$ tile $\{N+1,N+2,\dots\}$. A similar adjustment allows one to
easily change the range of $n$ in the definition of $\B$.

Skolem~\cite{MR19:1159i} stated (incorrectly) that if $\al,\be$ are positive irrationals, then
$\B(\al,\alpr)$ and $\B(\be,\bepr)$ tile $\{n \in\Z \colon n\ge
\min\{\floor{\frac{1-\alpr}{\al}},\floor{\frac{1-\bepr}{\be}}\}\}$ if and only if $\al+\be=1$ and
$\alpr+\bepr \in\Z$. Borwein \& Borwein~\cite{MR94i:11044} assume that $0<\al<1$, $\al$ irrational,
$0<\al+\alpr\le \al$, and $\frac{n-\alpr}{\al}$ is never integral. Under these hypotheses, they prove
(correctly) that $\B(\al,\alpr)$ and $\B(\be,\bepr)$ tile $\N$ if and only if $\al+\alpr=1$ and
$\alpr+\bepr=0$.

Fraenkel's statement (using a {\em simplified}\ form of his notation) is as follows. Let $\al$ and $\be$ be
positive real numbers, either both rational or both irrational, and let $\gamma$ and $\delta$ be arbitrary
real numbers. Let $S$ and $T$ be the sets of integers of the form $\phi_n=[n\al+\gamma]$ and
$\psi_n=[n\be+\delta]$, respectively, where $n$ ranges over $\N$. Further, assume that $\phi_1\le \psi_1$. If
$\al,\be$ are irrational, then $S$ and $T$ tile $\{n : n\ge \phi_1\}$ if and only if $\frac 1\al + \frac 1\be
=1$, $\frac{\gamma}{\alpha}+\frac{\delta}{\beta}=\phi_1-1$, and
    $$n\beta+\delta=K,\quad \text{$n,K$ integral implies $n<1$}.$$
If $\al,\be$ are rational, then $S$ and $T$ tile $\{n : n\ge \phi_1\}$ if and only if
$\frac{1}{\alpha}+\frac{1}{\beta}=1$ and
    $$\frac{\gamma}{\alpha}+\frac{\delta}{\beta}=\phi_1-1-a^{-1}+\eta+\rho,$$
where $\alpha=\frac ac$, $\gcd(a,c)=1$, $\frac{\gamma}{\alpha}\equiv \eta \pmod{a^{-1}}$, $0 \le
\eta<a^{-1}$, $\beta=\frac bd$, $\gcd(b,d)=1$, $\frac{\delta}{\beta}\equiv \rho \pmod{b^{-1}}$, $0\le \rho <
b^{-1}$. It is remarkable how much simplification we purchase by considering $\floor{\tfrac{n-\alpr}{\al}}$
in place of $\floor{n\al+\gamma}$. Our less-obvious definition leads to a somewhat simpler statement, and a
much simpler proof. (Note: our additional restriction to $0\le\al+\alpr\le 1$ in the irrational case and
$\frac 1q \le \al+\alpr\le 1$ in the rational case correspond to insisting that $\phi_1=1$.)

\section{Statement of Brown's Decomposition}\label{Sec.Statement.Brown}

Before stating Brown's Decomposition, we must define the {\em continued fraction expansion} of a natural
number. Let $[0;a_1,a_2,\ldots]$ be the continued fraction expansion of $\al$, and denote the continuants
(the denominators of the convergents to $\al$) by $q_0,q_1,\ldots$, i.e., $q_0:=1$, $q_1:=a_1$, and $q_i:=a_i
q_{i-1}+q_{i-2}$. For a positive integer $m$, define $z_0, z_1,\dots$ by writing $m$ greedily as a sum of
$q_i$ (that is, always use the largest $q_i$ possible):
    $$m=z_t q_t+z_{t-1}q_{t-1} + \dots z_1 q_1+z_0q_0.$$
We call $(z_tz_{t-1}\cdots z_1z_0)_\al$ the continued fraction expansion of $m$ with respect to $\al$. The
standard reference for this and other systems of numeration is \cite{MR86d:11016}.

We can now state Brown's Decomposition.

\begin{namedtheorem}{Brown's Decomposition}
Let $\al=[0;a_1,a_2,\dots]$ have continuants $q_0,q_1,q_2,\dots,$ and let $m=(z_t\cdots z_1z_0)_{\al}$. Then
for $i\ge 2$
    $$C_{q_i}=C_{q_{i-1}}^{a_i}C_{q_{i-2}} \quad \text{and} \quad
    C_m = C_{q_t}^{z_t} C_{q_{t-1}}^{z_{t-1}} \cdots C_{q_1}^{z_1}C_{q_0}^{z_0}.$$
\end{namedtheorem}

An avid reader may enjoy proving that $z_0,z_1,\dots$ is the unique sequence of nonnegative integers such
that $m=\sum_{i\geq0} z_i q_i$, $0\leq z_i \leq a_i$, and $(z_i=a_i) \Rightarrow (i\geq 1 \text{ and }
z_{i-1}=0)$. We remark that this is sometimes referred to as the Ostrowski expansion, and sometimes as the
Zeckendorf expansion, especially when $q_0, q_1, \dots$ are Fibonacci numbers.

\section{Preliminaries}\label{Sec.Preliminaries}

Much of our work will take place in $\R/\Z$. While there is no natural linear ordering of $\R/\Z$, there
is a natural `ternary' order: we say that real numbers $w, x, y$ are {\em in order} if there is a
nondecreasing function $f:[0,1]\to\R$ with $\fp{f(0)}=\fp{w}$, $\fp{f(\frac12)}=\fp{x}$,
$\fp{f(1)}=\fp{y}$, and $f(1)-f(0)<1$.

This definition is precise but awkward; there is a geometric description that is conceptually simpler. Let
$\tau:\R\to\C$ be defined by $\tau(z)=e^{2\pi i z}$. The range of $\tau$ is the circle $D:=\{ z \colon
|z|=1\}$, and the group $(D,\cdot)$ is isomorphic to $\R/\Z$ (in fact, $\tau$ is one isomorphism). We say
that $w,x,y$ are in order if $\tau(x)$ is on the counter-clockwise arc from $\tau(w)$ to $\tau(y)$. We write
$x\equiv y$ when $\tau(x)=\tau(y)$, i.e., when $x-y\in\Z$.

We define the arcs $\arc{(w,y)}, \arc{(w,y]}$, and $\arc{[w,y)}$ to be $\emptyset$ if $w\equiv y$, and
otherwise define the arcs through
    \begin{align*}
    \arc{(w,y)}     &:= \left\{ x \in \R \colon x\not\equiv w, x\not\equiv y,
                                \text{ and $w,x,y$ are in order.}\right\} \\
    \arc{(w,y]}     &:= \left\{ x \in \R \colon x\not\equiv w,
                                \text{ and $w,x,y$ are in order.}\right\} \\
    \arc{[w,y)}     &:= \left\{ x \in \R \colon x\not\equiv y,
                                \text{ and $w,x,y$ are in order.}\right\}
    \end{align*}

Our proofs of both Fraenkel's Theorem and Brown's Decomposition rely heavily on the following lemma.
\begin{lem}\label{BeattyTest}
Let $k$ be an integer, and $0<\al<1$. Then $k\in \B(\al,\alpr)$ if and only if
    $$\left( k\al+\alpr>0 \right)
        \quad \text{\tt AND}\quad
        \left( k\al \in\arc{(-\al-\alpr,-\alpr]}\right).$$
\end{lem}

\begin{proof}
\begin{align*}
    k\in \B(\al,\alpr)
        & \Leftrightarrow
            \exists n\in\N \left(k\le \tfrac{n-\alpr}{\al} < k+1 \right) \\
        & \Leftrightarrow
            \exists n\in\N \left(k\al+\alpr \le n < (k+1)\al+\alpr \right)\\
        & \Leftrightarrow
            \left(k\al+\alpr>0 \right) \quad\text{\tt AND}\quad \left((k+1)\al+\alpr\in
                \arc{(0,\al]}\right) \\
        & \Leftrightarrow
           \left(k\al+\alpr>0 \right) \quad\text{\tt AND}\quad
           \left(  k\al \in\arc{(-\al-\alpr,-\alpr]}\right)
\end{align*}
\end{proof}

\section{Proof of Fraenkel's Theorem}\label{Sec.Proof.Fraenkel}

\subsection{Spirit of Proof}\label{Sec.Special.Fraenkel}

In this section we prove a theorem whose statement and proof are similar to Fraenkel's Theorem, but for which
the technical details are considerably reduced. We note that Fraenkel~\cite{Fraenkel} also proved this
result. Set
    $$
    S_\al:=\left(\floor{\tfrac{n-\alpr}{\al}}\right)_{n=-\infty}^\infty
    \quad \text{and} \quad
    S_\be:=\left(\floor{\tfrac{n-\bepr}{\be}}\right)_{n=-\infty}^\infty.
    $$

\begin{thm}
Let $\al,\be$ be positive irrationals. The sequences $S_\al$ and $S_\be$ tile $\Z$ if and only if
$\al+\be=1$, $\alpr+\bepr\in\Z$, and $k\al+\alpr$ is never an integer (for $k\in\Z$).
\end{thm}

The additional difficulty of proving Fraenkel's Partition is in dealing with the `edge effects' introduced by
restricting the index $n$ in the definition of $S_\al$ and $S_\be$ to $n\ge1$, and in dealing with rational
$\al$. The proof of Fraenkel's Theorem given below is self-contained; this subsection is included only to
give the look-and-feel of our approach.

\begin{proof}
First, we note that the density of $S_\al$ is $\al$, and that of $S_\be$ is $\be$; it is clearly
necessary that $\al+\be=1$. From this point forward, we assume that $\al+\be=1$.

Next, observe that $k\in S_\al$ exactly if there is an $n\in\Z$ with $k \leq \frac{n-\alpr}{\al} < k+1$,
which is the same as
    $$k\al+\alpr \le n < k\al+\alpr+\al.$$
This, in turn, is the same as ${k\al+\alpr} \in \arc{(-\al,0]}.$ Thus (arguing identically for $\be$),
    $$S_\al = \left\{k \colon  {k\al}  \in \arc{(-\al-\alpr,-\alpr]}\right\}
        \text{\quad and \quad}
      S_\be = \left\{k \colon  {k\be}  \in \arc{(-\be-\bepr,-\bepr]}\right\}.$$
But $ {k\be} =  {k(1-\al)} \equiv - {k\al} $, so that $S_\be$ is also given by
    $$S_\be = \left\{k \colon  {k\al} \in \arc{[\bepr,\be+\bepr)} \right\}.$$
Thus, $k\in S_\al$ if ${k\al}  \in \arc{(-\al-\alpr,-\alpr]}=:{A}$ and $k\in S_\be$ if ${k\al} \in
\arc{[\bepr,\be+\bepr)}=:{B}$. Since $(\fp{k\al})_{k=-\infty}^\infty$ is dense in $[0,1)$, if ${A}$ and ${B}$
intersect in an arc, there are infinitely many $k$ in both $S_\al$ and $S_\be$; and if ${A}$ and ${B}$ both
omit some arc, there are infinitely many $k$ in neither $S_\al$ nor $S_\be$. It follows that the right
endpoint of $A$ is the left endpoint of ${B}$, i.e., $-\alpr\equiv \bepr$, which is the same as $\alpr+\bepr
\in \Z$. The only point in $A\cap B$ is $-\alpr\equiv \be$, so that if ${k\al} \equiv -\alpr\equiv \be$, then
$k$ is in both $S_\al$ and $S_\be$. This happens if and only if $k\al+\alpr$ is an integer.
\end{proof}

\subsection{Without loss of generality ...}

\subsubsection{Fraenkel's Partition is symmetric in $\al$ and $\be$.}\label{Sec.sub.sub.symmetry}
We first note that the Theorem is symmetric in $\al$ and $\be$. Obviously, ``The sequences $\B(\al,\alpr)$
and $\B(\be,\bepr)$ tile $\N$'' has the claimed symmetry. Combining Conditions~{\em 2} and~{\em 3}, we find
$0<\be<1$, the symmetric counterpart to Condition~{\em 1}. Condition~{\em 2} is symmetric as stated. We
return to Condition~{\em 3} in the next paragraph. In Condition~{\em 4}, the equation $\alpr+\bepr=0$ is
already symmetric; and ``$k\al+\alpr \not \in \Z$'' implies that
$k\be+\bepr=k(1-\al)+(-\al)=k-(k\al+\alpr)\not\in\Z$ (using Condition~\emph{2} and $\alpr+\bepr=0$). Thus
Condition~\emph{4} implies its symmetric counterpart. If $\al$ is rational, then $\be=1-\al$ (from
Condition~\emph{2}) is also, and moreover both $\al$ and $\be$ have the same denominator. Thus
Condition~\emph{5} implies its symmetric counterpart also.

If $\al$ is irrational, then $\al=1-\be$ (from Condition~\emph{2}) and $\alpr=-\bepr$ (from
Condition~\emph{4}), so $0\le \al+\alpr \le 1$ (from Condition~\emph{3}) implies $0\le \be+\bepr \le 1$, the
symmetric twin of Condition~\emph{3}. Now suppose that $\al$ is rational, say $q\al\in\N$. Then the
inequalities $\tfrac 1q \le \al+\alpr \le 1$ (from Conditions~\emph{3} and~\emph{5}) yield $1\le
q\al+q\alpr\le q$, and since $1,q\al,q$ are all integers, this yields $1\le q\al+\ceiling{q\alpr} \le q$. Now
from Condition~\emph{2}, $\al=1-\be$, and from Condition~\emph{5}, $\ceiling{q\alpr}=1-\ceiling{q\bepr}$, so
the inequalities imply $1\le q(1-\be)+1-\ceiling{q\bepr}\le q$, which simplifies to $1\le
q\be+\ceiling{q\bepr}\le q$. Now, since $q\be=q-q\al\in\N$, and $1,q$ are integers also, the inequalities
become $1\le q\be+q\bepr \le q$, or more simple $\tfrac 1q \le \be+\bepr \le 1$. This gives the symmetric
counterpart to Condition~\emph{3} (since $\frac 1q >0$) and to Condition~\emph{5}.

\subsubsection{If $\al$ is rational, then $\al,\alpr,\be,\bepr$ are all rational with the same
denominator.}\label{Sec.sub.sub.rational.alpha}

We now observe that if $\al$ is rational (with denominator $q$), we can assume without loss of generality
that $\alpr$ is also rational with denominator $q$. For $0\le \al+\alpr \le 1$ if and only if $0\le
\al+\alprpr \le 1$ and clearly $\ceiling{q \alprpr}+\ceiling{q \alprpr}=\ceiling{q\alpr}+\ceiling{q\bepr}$,
and so the Conditions~\emph{1--5} are unaffected by replacing $\alpr$ with $\alprpr$. By
Lemma~\ref{simplify.rationals} below, the sequence $\B(\al,\alpr)$ is also unaffected. Thus, we assume from
this point on that if $\al$ is rational with denominator $q$, then so is $\alpr$. Further, if $\be$ is
rational with denominator $q$, we assume that $\bepr$ is too. The equation in Condition~\emph{5} now
simplifies to $\alpr+\bepr=\frac 1q$.

When $\al,\be$ are known to be positive rationals with $\al+\be=1$, we define natural numbers $q,a,\ap,b,\bp$
by the conditions
    $$
    \al=\frac aq, \quad
    \gcd(a,q)=1, \quad
    \alpr = \frac{\ap}{q}, \quad
    \be = \frac bq, \quad
    \bepr = \frac{\bp}{q}.
    $$

\begin{lem}\label{simplify.rationals}
For any $a,q\in \N$ and $\alpr\in \R$, $\B(\frac aq, \alpr)=\B(\frac aq, \frac{\ceiling{q\alpr}}{q})$.
\end{lem}

\begin{proof}
Set $\al=\frac aq$. We have $k\al+\alpr>0 \Leftrightarrow ka+q\alpr>0$, and since $ka\in\Z$, we have
$ka+q\alpr>0 \Leftrightarrow ka+\ceiling{q\alpr}>0 \Leftrightarrow k\al + \alprpr>0$.

Now $k\al$ is rational with denominator $q$, and there are no such numbers in $\arc{(-\alprpr,-\alpr]}$ or in
$\arc{(-\al-\alprpr,-\al-\alpr]}$. Thus
    $$k \al \in \arc{(-\al,-\al-\alpr]} \Leftrightarrow
        k\al \in \arc{(-\al-\tfrac{\ceiling{q\alpr}}{q},-\tfrac{\ceiling{q\alpr}}{q}]}.$$
Apply Lemma~\ref{BeattyTest} to finish the proof.
\end{proof}

\subsubsection{The sets $A$ and $B$ are important.}
We define
    $$ A:= \arc{(-\al-\alpr,-\alpr]}
        \text{ and }
       B:= \arc{[\bepr,\be+\bepr)}.
    $$
\begin{lem}\label{lemma.AandB}
Suppose $\al+\be=1$. Then $k\in \B(\al,\alpr)$ if and only if
    $$
        \left( k\al+\alpr>0 \right)
        \quad \text{\tt AND}\quad
        \left( k\al \in A\right).
    $$
Further, $k\in \B(\be,\bepr)$ if and only if
    $$
        \left( k\be+\bepr>0 \right)
        \quad \text{\tt AND}\quad
        \left( k\al \in B\right).
    $$
\end{lem}

\begin{figure}
    \begin{center}
        \begin{picture}(225,225)
            \put(35,135){$A$}
            \put(175,105){$B$}
            \ifpdf
                \put(0,-70){\includegraphics{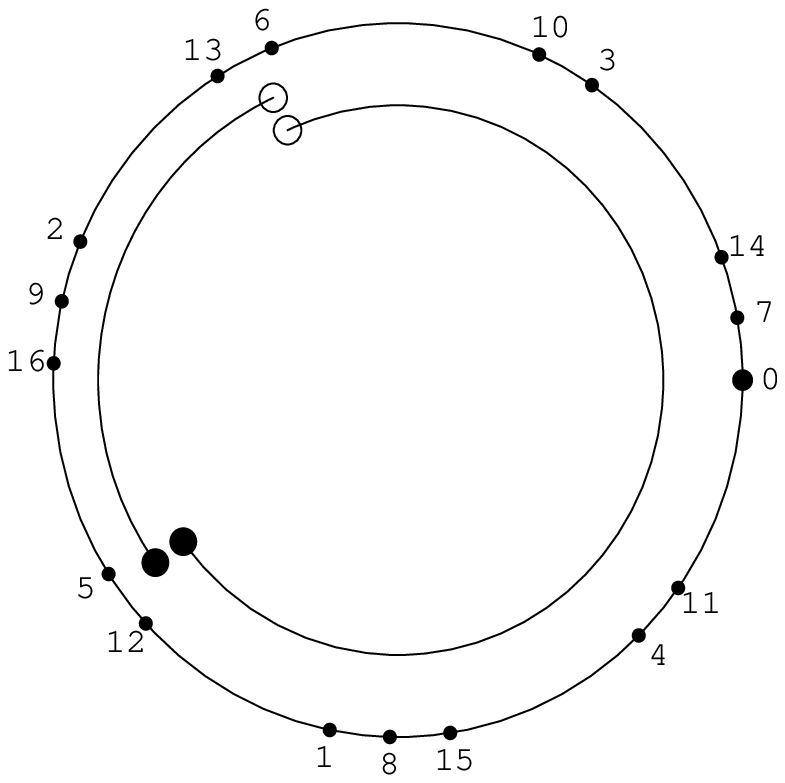}}
            \else
                \put(0,0){\includegraphics{PartitionPlot}}
            \fi
        \end{picture}
    \end{center}
    \caption{The partition of $\N$ into $\B(3-e,\tfrac25)=\{2,5,9,13,16,\dots\}$ and
    $\B(e-2,-\tfrac 25)=\{1,3,4,6,7,8,\dots\}$.}
    \label{figure.partition}
\end{figure}

Figure~\ref{figure.partition} gives an example of Fraenkel's Partition with $\al=e-2$, $\alpr=\tfrac23$,
$\be=3-e$ and $\bepr=-\tfrac23$. Shown in the figure are the circle $|z|=1$ and the angles that correspond to
the sets $A$ and $B$. Also shown are the points $\tau(k\al)$ ($0\le k \le 16$), labelled ``$k$''.

\begin{proof}
The condition for $k\in \B(\al,\alpr)$ is simply a restatement of Lemma~\ref{BeattyTest}. To prove the
condition for $k\in\B(\be,\bepr)$, we must show that
    $$k\be \in \arc{(-\be-\bepr,-\bepr]}
    \Leftrightarrow
    k\al \in \arc{[\bepr,\be+\bepr)}.
    $$
This is obvious since $k\be=k(1-\al) \equiv -k\al$, and from the observation that $-w,-k\al,-y$ are in order
if and only if $y, k\al, w$ are too.
\end{proof}

\subsection{Conditions~{\em 1--5} are Sufficient}

Suppose that $\al,\alpr,\be,\bepr$ satisfy Conditions~{\em 1--5}. We have $\frac{n-\alpr}{\al}\ge
\frac{1-\alpr}{\al}\ge1$ (using Condition~\emph{3}), so $B(\al,\alpr)\subseteq \N$, and likewise (by the
symmetry proved in Section~\ref{Sec.sub.sub.symmetry}) $B(\be,\bepr)\subseteq\N$.

We now show the equivalence (for $k\in\N$)
    $$k\in \B(\al,\alpr) \Leftrightarrow k\not\in \B(\be,\bepr),$$
i.e., we show that the two sequences tile $\N$. We break the work into four cases: $k=1$ or $k>1$, and
$\al\in\Q$ or $\al\not\in\Q$.

\subsubsection{$k=1$.} Since $\al\in(0,1)$ and $\be=1-\al\in(0,1)$, the two sequences are strictly
increasing. Since $B(\al,\alpr)$ and $B(\be,\bepr)$ are in $\N$ (as proved above), we have for
$x\in\{\al,\be\}$
    $$1\in B(x,x^\prime) \Leftrightarrow \frac{1-x^\prime}{x} < 2.$$
Specifically, we have
    \begin{equation}\label{k1Ba}
    1 \in \B(\al,\alpr) \Leftrightarrow \frac{1-\alpr}{\al}<2 \Leftrightarrow 1<2\al+\alpr
    \end{equation}
and
    \begin{equation}\label{k1Bb}
    1 \not \in \B(\be,\bepr) \Leftrightarrow \frac{1-\bepr}{\be} \ge 2
        \Leftrightarrow 1\ge 2\be+\bepr
    \end{equation}

If $\al$ is irrational then, by Condition~\emph{4}, $2\al+\alpr\not=1$ and so $1<2\al+\alpr \Leftrightarrow
1\le 2\al+\alpr$. By Condition~\emph{2}, $\al=1-\be$, and by Condition~\emph{4}, $\alpr=-\bepr$. Thus,
    $$
    1< 2\al+\alpr \Leftrightarrow 1\le 2\al+\alpr \Leftrightarrow 1\ge 2\be+\bepr,
    $$
connecting the equivalences in Lines~(\ref{k1Ba}) and~(\ref{k1Bb}).

Now suppose that $\al$ is rational, and recall that we may assume that $\al,\alpr,\be,\bepr$ are all
rationals with denominator $q$, as per the discussion in Section~\ref{Sec.sub.sub.rational.alpha}. Using
$\al=1-\be$ and $\alpr=\frac1q-\bepr$, we have
    $$1<2\al+\alpr \Leftrightarrow 2\be+\bepr<1+\frac1q \Leftrightarrow 2\be+\bepr \le 1,$$
the last equivalence following from $\be,\bepr$ being rationals with denominator $q$. This connects the
equivalences in Lines~(\ref{k1Ba}) and~(\ref{k1Bb}).

\subsubsection{$k>1$.}
Suppose that $k>1$, so that $k\al+\alpr>\al+\alpr\ge 0$ (by Condition~\emph{3}), and by symmetry
$k\be+\bepr>0$.

Lemma~\ref{lemma.AandB} reduces to
    $$
    k\in \B(\al,\alpr) \Leftrightarrow {k\al} \in A
    \quad \text{and} \quad
    k\in \B(\be,\bepr) \Leftrightarrow k\al \in B.
    $$
Thus, to prove $k\in \B(\al,\alpr) \Leftrightarrow k\not\in \B(\be,\bepr)$ it will suffice to show that $k\al
\not\in A\cap B$ and $k\al \not\in A^c \cap B^c$.

\begin{figure}
    \begin{center}
    \begin{picture}(215,180)
        \put(180,87){$0\equiv 1$}
        \put(40,40){$A$}
        \put(124,124){$B$}
        \put(-15,161){$-\al-\alpr\equiv $}
        \put(-10,149){$\be+\bepr$}
        \put(159,25){$-\alpr \equiv \bepr$}
    \ifpdf
        \put(0,0){\includegraphics{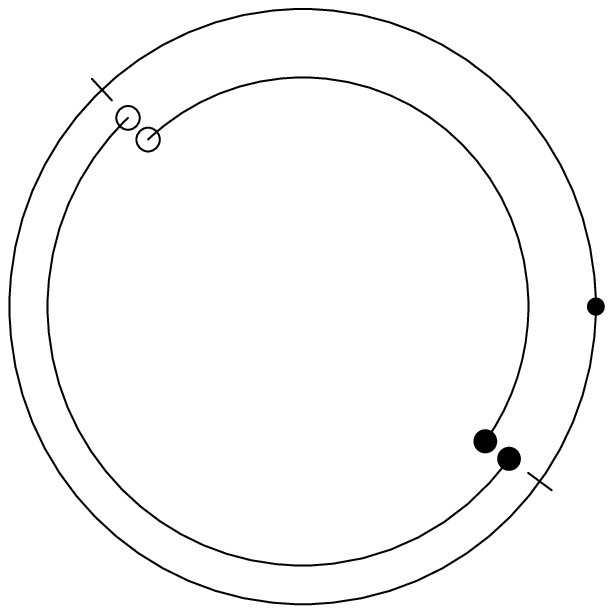}}
    \else
        \put(0,0){\includegraphics{IrrationalPlot}}
    \fi
    \end{picture}
    \caption{$A$ and $B$, with $\al$ irrational}
    \label{IrrationalCase}
    \end{center}
\end{figure}

\begin{figure}
    \begin{center}
    \begin{picture}(215,180)
        \put(180,87){$0\equiv 1$}
        \put(42,42){$A$}
        \put(122,122){$B$}
        \put(0,163){$-\al-\alpr$}
        \put(-15,135){$\be+\bepr$}
        \put(152,24){$-\alpr$}
        \put(173,56){$\bepr$}
    \ifpdf
        \put(0,0){\includegraphics{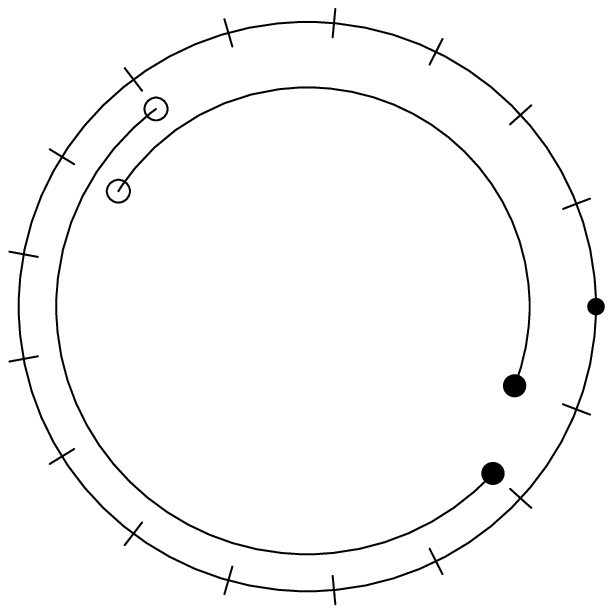}}
    \else
        \put(0,0){\includegraphics{RationalPlot}}
    \fi
    \end{picture}
    \caption{$A$ and $B$, with $\al$ rational}
    \label{RationalCase}
    \end{center}
\end{figure}

The next two paragraphs are accompanied by Figures~\ref{IrrationalCase} and~\ref{RationalCase}. In these
figures the point $\tau(z)$ is labelled ``$z$''. These figures show the circle $\tau(\R)$ (the outer circle
with the point $\tau(0)=\tau(1)$), and display the angles corresponding to $\tau(A)$ and $\tau(B)$ as arcs
inside the circle. Also labelled are the points $\tau(-\alpr), \tau(\bepr), \tau(-\al-\alpr),$ and
$\tau(\be+\bepr)$.

If $\al$ is irrational (see Figure~\ref{IrrationalCase}), then
$B:=\arc{[\bepr,\be+\bepr)}=\arc{[-\alpr,-\al-\alpr)}$. We have $A\cap B=\{y \colon y\equiv -\alpr\}$. Since
$k\al+\alpr \not\equiv 0$ by Condition~\emph{4}, we know that $k\al\not\in A \cap B$. We have $A^c \cap
B^c=\{y \colon y\equiv -\al-\alpr\}$. If $k\al \in A^c \cap B^c$ then $k\al \equiv -\al-\alpr$. Thus
$(k+1)\al+\alpr \in \Z$, but this is forbidden by Condition~\emph{4}.

If $\al$ is rational (see Figure~\ref{RationalCase}), then $B:=\arc{[\bepr,\be+\bepr)}=\arc{[\frac
1q-\alpr,\frac1q-\al-\alpr)}$. We have $A\cap
B=\arc{(-\al-\alpr,\be+\bepr)}=\arc{(\frac{-a-\ap}{q},\frac{-a-\ap+1}{q})}$ and $A^c \cap B^c =
\arc{(-\alpr,\bepr)}=\arc{(\frac{-a}{q},\frac{-a+1}{q})}$, neither of which contain a rational with
denominator $q$. Since $k\al$ is rational with denominator $q$, we know that $k\al \not \in A \cap B$ and
$k\al \not \in A^c \cap B^c$.

\subsection{Conditions~{\em 1--5} are Necessary}

We now assume that $\B(\al,\alpr)$ and $\B(\be,\bepr)$ tile $\N$, and prove Conditions~{\em 1--5}.

We define $k_0:= \max\{0,-\alpr/\al,-\bepr/\be\}$. For $k>k_0$, then both $k\al+\alpr$ and $k\be+\bepr$ are
positive, and this will simplify our work tremendously.

The sequences $\B(\al,\alpr), \B(\be,\bepr)$ contain infinitely many positive integers, so $\al$ and
$\be$ are nonnegative.

We can count the number of elements of $\B(\al,\alpr)$ which are less than an integer $k$:
 \begin{align*}
  \left| \B(\al,\alpr) \cap (-\infty ,k)\right| &=
    \left| \left\{ n\ge 1\colon \floor{\frac{n-\alpr}{\al}}<k\right\}\right| \\
    &= \left| \left\{ n\ge 1\colon \frac{n-\alpr}{\al}<k\right\}\right| \\
    &= \max\left\{ 0, \ceiling{k\al+\alpr}-1 \right\}
 \end{align*}
and likewise
 $$\left| \B(\be,\bepr) \cap (-\infty,k)\right|  = \max\left\{ 0, \ceiling{k\be+\bepr}-1 \right\}.$$
Since $\B(\al,\alpr)$ and $\B(\be,\bepr)$ tile $\N$, we have (for $k>k_0$)
 \begin{align}
  k-1
    &=  \left| (\B(\al,\alpr) \cup \B(\be,\bepr)) \cap (-\infty,k) \right| \notag \\
    &=  \left| \B(\al,\alpr) \cap (-\infty,k)\right| +\left| \B(\be,\bepr) \cap (-\infty,k)\right| \notag\\
    &=  \ceiling{k\al+\alpr}+\ceiling{k\be+\bepr}-2. \label{countingLine}
  \end{align}
Divide by $k$, and let $k$ go to infinity to find $1=\al+\be$ (Condition~\emph{2}).

As $\al$ and $\be$ are positive, and $\al+\be=1$, we find $0<\al<1$ (Condition~\emph{1}).

We will use the following short lemma several times.
\begin{lem}\label{aprimebprime}
If $k>k_0$, $\fp{k\al+\alpr}\not=0$, and $\fp{k\be+\bepr}\not=0$, then
    $$\fp{k\al+\alpr}+\fp{k\be+\bepr}=\alpr+\bepr+1.$$
\end{lem}

\begin{proof}
Since $k\al+\alpr$ is not an integer, we have $\ceiling{k\al+\alpr}=k\al+\alpr+1-\fp{k\al+\alpr}$, and
likewise $\ceiling{k\be+\bepr}=k\be+\bepr+1-\fp{k\be+\bepr}$. Equation~(\ref{countingLine}) now simplifies to
finish the proof.
\end{proof}

The hypothesis that $k\in \B(\al,\alpr) \Leftrightarrow k\not\in\B(\be,\bepr)$ and Lemma~\ref{AandNotB} imply
that (for $k>k_0$)
    $${k\al} \in A \Leftrightarrow {k\al} \not\in B.$$
This, in turn, is equivalent to the assertion for $k>k_0$
    \begin{equation}\label{AandNotB}
    {k\al} \not \in A \cap B
    \quad \text{and} \quad
    {k\al} \not \in A^c \cap B^c.
    \end{equation}

\subsubsection{If $\al$ is irrational ...}

The set $\{ \fp{k\al} \colon k>k_0\}$ is dense in $[0,1)$; Line~(\ref{AandNotB}) implies that there is no
nontrivial interval contained in $A\cap B$, nor in $A^c\cap B^c$. Thus, the right endpoint of $B$ is
congruent (modulo 1) to the left endpoint of $A$, i.e., $\bepr \equiv -\alpr$.

Suppose (by way of contradiction) that both $s_1\al+\alpr$ and $s_2\al+\alpr$ are integers (with $s_i\in\Z$).
Then so is $(s_1-s_2)\al=(s_1\al+\alpr)-(s_2\al+\alpr)$, contrary to our hypothesis that $\al$ is irrational.
Thus, there is at most one integer $s$ such that $s\al+\alpr\in\Z$. Consequently, we may choose $k_1$ so that
neither $k_1\al+\alpr$ nor $k_1 \be+\bepr$ are integers. We know from the previous paragraph that
$\bepr+\alpr$ is an integer, so from Lemma~\ref{aprimebprime}, $\fp{k_1\al+\alpr}+\fp{k_1\be+\bepr}$ is an
integer. By the definition of fractional part, $0 \le \fp{k_1\al+\alpr}+\fp{k_1\be+\bepr}<2$, and by the
choice of $k_1$, $\fp{k_1\al+\alpr}+\fp{k_1\be+\bepr}>0$. Thus $\fp{k_1\al+\alpr}+\fp{k_1\be+\bepr}=1$, and
Equation~(\ref{AandNotB}) reduces to $\alpr+\bepr=0$.

Now suppose that $k\al+\alpr \in \Z$ (with $k\ge2$, not necessarily larger than $k_0$). Note that $k\be+\bepr
=k(1-\al)-\alpr =k-(k\al+\alpr)\in\Z$ also. We have $k\al\equiv -\alpr \equiv \bepr \in A\cap B$, so
Lemma~\ref{AandNotB} implies that $k\al+\alpr\le 0$ or $k\be+\bepr \le 0$. Without loss of generality,
suppose that $k\al+\alpr \le 0$. Then $(k-1)\al+\alpr\le 0$, so by Lemma~\ref{BeattyTest}, $k-1 \not \in
\B(\al,\alpr)$. Also $(k-1)\al \equiv -\al-\alpr \equiv \be+\bepr \not \in B$, so $k-1 \not\in
\B(\be,\bepr)$. That is, $k-1$ is in neither $\B(\al,\alpr)$ nor $\B(\be,\bepr)$, contradicting the
hypothesis that these sets tile $\N$. This establishes Condition~\emph{4}.

Since $\B(\al,\alpr)\subseteq\N$, we have $\frac{1-\alpr}{\al}\ge1$, whence $\al+\alpr\leq 1$. Likewise,
$1\ge \be+\bepr =1-\al-\alpr$, whence $0\le \al+\alpr$, establishing Condition~\emph{3}.

\subsubsection{If $\al$ is rational ...}

If $\al=\frac12$,then $\be=1-\al=\frac12$. We have $\B(\al,\alpr)=(2n-\ap)_{n=1}^\infty$ and
$\B(\be,\bepr)=(2n-\bp)_{n=1}^\infty$. This is the even-odd tiling of $\N$; we have $\{\ap,\bp\}=\{0,1\}$.
Conditions~\emph{3} and~\emph{5} are now easily verified. From this point on, we can assume that one of
$\al,\be$ is strictly less than $\frac12$.

The set of fractional parts $\{ \fp{k\al} \colon k>k_0\}= \{ \tfrac iq \colon 0\le i < q\}$, so
Line~(\ref{AandNotB}) says that there is no multiple of $\tfrac 1q$ in $A\cap B$ or in $A^c \cap B^c$.
Consider the multiples of $\tfrac 1q$
    $$-\al-\alpr, -\al-\alpr+\tfrac1q, -\al-\alpr+\tfrac 2q, \dots, -\alpr,-\alpr+\tfrac1q.$$
Since $\al \le 1-\tfrac1q$, the first and the last (which may be congruent modulo 1) are clearly not in
$A:=\arc{(-\al-\alpr,-\alpr]}$. Therefore, they must be in $B:=\arc{[\be+\bepr,\bepr)}$. That is
$B=\arc{[-\alpr+\tfrac1q,-\al-\alpr+\tfrac1q)}$. In particular, $\bepr \equiv -\alpr+\frac1q$.

If $x\not\in\Z$ then for any $k$ at most one of $s x+x^\prime$, $(s+1)x+x^\prime$ is integral (where
$s\in\Z$), since their difference is not integral. If $0<x<\frac12$, then at most one of $s x+x^\prime$,
$(s+1)x+x^\prime$, $(s+2)x+x^\prime$ is integral (if two were, then their difference would be also, but their
difference is $<1$). Thus, since one of $\al,\be$ is strictly less than $\frac12$, we may choose $k_1 \in
\{\floor{k_0+1},\floor{k_0+2},\floor{k_0+3}\}$ with neither $k_1\al+\alpr$ nor $k_1\be+\bepr$ integral. In
particular $\fp{k_1\al+\alpr}+\fp{k_1\be+\bepr}\ge \frac2q$.

Choose $k_1>k_0$ so that $\fp{k_1\al+\alpr}$ and $\fp{k_1\be+\bepr}$ are positive. By
Lemma~\ref{aprimebprime} we have $\fp{k_1\al+\alpr}+\fp{k_1\be+\bepr}-\frac1q \in \Z$, and by the definition
of fractional part, $\fp{k_1\al+\alpr}+\fp{k_1\be+\bepr}-\frac1q \in [-\frac1q,2-\frac1q)$. Thus,
$\fp{k_1\al+\alpr}+\fp{k_1\be+\bepr}-\frac1q =1$. Plugging this into Lemma~\ref{aprimebprime}, we find
$\alpr+\bepr=\frac1q$.

Since $\B(\al,\alpr)\subseteq \N$, we have $\tfrac{1-\alpr}{\al}\ge1$, whence $\al+\alpr \le 1$. Likewise,
$1\ge \be+\bepr=1-\al+\frac1q-\alpr$, so that $\al+\alpr \ge \frac1q\ge 0$. This establishes
Condition~\emph{3} and the last piece of Condition~\emph{5}.

\subsection{Using Nonstandard Analysis to Derive Irrational Case}\label{Sec.nonstandard}

We can use nonstandard analysis to easily derive the irrational case of Fraenkel's Partition from the
rational case. Specifically, we now show that if Conditions~\emph{1--3},~\emph{5} are sufficient for rational
$\al$, then Conditions~\emph{1--4} are sufficient for irrational $\al$.

Suppose that $\al$ is irrational and $\al,\alpr,\be,\bepr$ satisfy Conditions~\emph{1--4}. Since
Conditions~\emph{1--4} are symmetric in $\al,\be$ (as per the comment in Section~\ref{Sec.sub.sub.symmetry}),
we may label $\al$ and $\be$ so that $\al+\alpr\not=0$ and $\be+\bepr\not=1$.

Let $q\in {}^\ast\N \setminus \N$ be an infinite prime number. Set
    $$
    \gamma := \frac{\floor{q\al}}{q},
    \quad
    \gamma^\prime := \frac{\floor{q\alpr}}{q},
    \quad
    \delta := \frac{\ceiling{q\be}}{q},
    \quad
    \delta^\prime := \frac{\ceiling{q\be}+1}{q},
    $$
all $\ast$-rationals with denominator $q$. It is straightforward to verify that $\gamma,\gamma^\prime,
\delta, \delta^\prime$ satisfy Conditions~\emph{1}, \emph{2}, \emph{3}, and \emph{5}, so
$\B(\gamma,\gamma^\prime)$ and $\B(\delta,\delta^\prime)$ tile $^\ast\N$.

We have $\frac{n-\gamma^\prime}{\gamma}\ge \frac{n-\alpr}{\al}$ and for $n$ finite
$st(\frac{n-\gamma^\prime}{\gamma})=\frac{n-\alpr}{\al}$, so
$\floor{\frac{n-\gamma^\prime}{\gamma}}=\floor{\frac{n-\alpr}{\al}}$. Thus $B(\al,\alpr)=\N\cap
\B(\gamma,\gamma^\prime)$. We also have $st(\frac{n-\delta^\prime}{\delta})=\frac{n-\bepr}{\be}$, so
$\floor{\frac{n-\delta^\prime}{\delta}}=\floor{\frac{n-\bepr}{\be}}$ unless $\frac{n-\bepr}{\be}$ is an
integer. If $\frac{n-\bepr}{\be}=k \ge 2$, then $k-n=k-(k\be+\bepr)=k\al+\alpr \in\Z$, contrary to
Condition~\emph{4}. If $\frac{n-\bepr}{\be}=1$, then $\be+\bepr=1$, contrary to our assumption that
$\be+\bepr\not=1$.

\section{Proof of Brown's Decomposition}\label{Sec.Proof.Brown}

Fix an irrational $\al=[0;a_1,a_2,\dots]$, and let $q_0,q_1,\dots$ be its continuants. Specifically, let
$q_t$ be the largest continuant strictly less than $m$. We denote the distance to the nearest integer by
    $$\|x\| := \min\big\{x-\floor{x},\ceiling{x}-x\big\}=\min\big\{\fp{x},1-\fp{x}\big\}.$$

\begin{lem}\label{simple.Brown}
$C_m=C_{q_t}C_{m-q_t}.$
\end{lem}

Brown's Decomposition follows immediately from Lemma~\ref{simple.Brown} by induction. The proof of
Lemma~\ref{simple.Brown} relies on the following well-known result from the theory of continued fractions. It
can be found as Theorem 10.15 (page 370) of \cite{MR89b:11002}, for example, or on page 163 of
\cite{MR33:3981}.

\begin{lem}\label{spacing}
If $|s|< q_{t+1}$ then $\|s\al\| > \|q_t\al\|$.
\end{lem}

\begin{proofof}{Lemma~\ref{simple.Brown}}
First, if $\al=[0;a_1,\dots,a_n]$ is rational, then we can replace it with an irrational $x$ between $\al$
and $\al^\prime=[0;a_1,\dots,a_n,m+1]$. The continuants of $\al^\prime$ that are less than $m$ are the same
as the continuants of $\al$, so the continued fraction expansion of $m$ looks the same. And if the irrational
$x$ is sufficiently close to $\al$, then $\B(\al)\cap[1,m]=\B(x)\cap[1,m]$ and in particular $C_m$ does not
change when $\al$ is replaced with $x$. Thus, we may assume that $\al$ is irrational.

If $q_t=1$, then $t<2$ since $q_2=a_2q_1+q_0\ge2>q_t$. Thus either $m=(m)_\al$ or $m=(m0)_\al$. In the first
case ($t=0$), we have $m<q_1=a_1$, and so $\al<\frac 1{a_1}<\frac 1m$. Therefore $C_m=0^m=C_0 C_{m-1}$,
proving the Lemma. In the second case ($t=1$), we have $a_1=1$ and $1\le m <q_2=a_2+1$, so
$\al>\frac{a_2}{a_2+1}$. Now ${k}/{\al}<k \frac{a_2+1}{a_2}=k+\frac {k}{a_2}$, so $\floor{k/\al}=k$ for
$k=m\le a_2$. Thus $C_m=1^m=C_1C_{m-1}$, proving the Lemma. Thus we may assume from this point on that
$q_t>1$.

It is obvious that $C_{q_t}$ is an initial segment of $C_m$; the content of Lemma~\ref{simple.Brown} is that
$c_{q_t+k}=c_k$ for $1\le k \le m-q_t\le q_{t+1}-q_t$. In particular both $|k|$ and $|k+1|$ are strictly less
than $q_{t+1}$.

By definition
    $$c_k=1  \quad \Leftrightarrow \quad k\in \B(\al)$$
and by Lemma~\ref{BeattyTest} (with $\alpr=0$ and $k\ge1$)
    $$k\in \B(\al) \quad \Leftrightarrow \quad {k\al} \in \arc{(-\al,0]}.$$
Likewise,
    \begin{equation*}
    c_{q_t+k}=1
        \Leftrightarrow q_t+k \in \B(\al)
        \Leftrightarrow (q_t+k)\al \in \arc{(-\al,0]}
        \Leftrightarrow {k\al} \in \arc{(-\al-q_t\al,-q_t\al]}.
    \end{equation*}
We need to show only that $k\al$ is not in the symmetric difference of $\arc{(-\al,0]}$ and
$\arc{(-\al-q_t\al,-q_t\al]}.$ This symmetric difference is contained in
    $$\arc{[-\|q_t\al\|,\|q_t\al\|]} \cup \arc{[-\|q_t\al\|-\al,\|q_t\al\|-\al]}.$$
Now $|k|<q_{t+1}$, so by Lemma~\ref{spacing}, $\|k\al\|>\|q_t\al\|$ and consequently
    $$k\al \not \in \arc{[-\|q_t\al\|,\|q_t\al\|]}.$$
Also $|k+1|<q_{t+1}$, so $\|(k+1)\al\|>\|q_t \al\|$, which is the same as $(k+1)\al \not\in
\arc{[-\|q_t\al\|,\|q_t\al\|]}$. Thus
    $$k\al \not \in  \arc{[-\|q_t\al\|-\al,\|q_t\al\|-\al]}.$$
\end{proofof}

\end{document}